\documentclass[12pt]{article}
\usepackage{amssymb}
\usepackage{amsmath}

\newtheorem{theorem}{Th\'{e}or\`{e}me}
\newtheorem{corollary}[theorem]{Corollaire}
\newtheorem{definition}[theorem]{D\'{e}finition}
\newtheorem{lemma}[theorem]{Lemme}
\newtheorem{proposition}[theorem]{Proposition}

\def\Rea#1{\mathop{\rm Re}\nolimits}
\def\R{\mathbb{R}}
\def\FdeA{\mathcal{F}(\overline{A})}
\def\Fde#1{\mathcal{F}(#1)}
\def\GdeA{\mathcal{G}(\overline{A})}
\def\Gde#1{\mathcal{G}(#1)}

\begin{document}

\title{Une remarque sur les espaces d'interpolation faiblement localement
uniform\'{e}ment convexes}
\author{Daher Mohammad \\
D\'{e}partement de Math\'{e}matiques, Universit\'{e} Paris 7 \\
e-mail: m.daher@orange.fr}
\maketitle

\noindent\ \ \ \ \ \ \emph{R\'{e}sum\'{e}}: Soient $(A_0, A_1)$ un couple
d'interpolation, et $B_j$ l'adh\'{e}rence de $A_0^\ast \cap A_1^\ast$ dans
$A_j^\ast$, $j = 0, 1$. Pour tout $\theta \in \, ]0, 1[$, il existe une
contraction injective naturelle $R^\theta :
 A^\theta \rightarrow (B_0^\ast, B_1^\ast)^\theta$. On suppose que, pour
un $\beta \in \, ]0, 1[$, l'adh\'{e}rence de $R^\beta (A^\beta)$ dans 
$(B_0^\ast, B_1^\ast)^\beta$ est faiblement LUR. Alors 
$A^\theta = A_\theta$ pour tout $\theta \in \, ]0, 1[$.

\emph{Abstract}: Let $(A_0, A_1)$ be an interpolation couple, and let $B_j$
be the closure of $A_0^\ast \cap A_1^\ast$ in $A_j^\ast$, $j = 0, 1$. For
every $\theta \in \, ]0, 1[$, there exists a natural one to one contraction
$R^\theta : A^\theta \rightarrow (B_0^\ast, B_1^\ast)^\theta$. For some
$\beta \in \, ]0, 1[$, the closure of $R^\beta (A^\beta)$ in
$(B_0^\ast, B_1^\ast)^\beta$ is supposed to be weakly LUR. Then 
$A^\theta = A_\theta$ for every $\theta \in \, ]0, 1[$.


\emph{Mots cl\'{e}s}: Interpolation, espace faiblement-LUR

\emph{Avertissement}: Une premi\`{e}re version de ce travail a \'{e}t\'{e}
publi\'{e}e dans Colloq. Math. Vol. 113, No. 2, 197-204, (2011). Le lemme~1
de cette version est malheureusement faux, ce qui oblige \`{a} corriger
l'ensemble. L'auteur remercie Sten Kaijser de lui avoir signal\'{e} l'erreur
dans la preuve. Un rectificatif indiquant les corrections a \'{e}t\'{e}
envoy\'{e} \`{a} Colloq. Math.. Il nous a cependant sembl\'{e} utile de
pr\'{e}senter une version r\'{e}vis\'{e}e compl\`{e}te.

\begin{center}
\medskip
\end{center}

\section{Introduction et notations}

On note $X^\ast$ le dual d'un espace de Banach $X$.

Soit $\overline{A} = (A_0, A_1)$ un couple d'interpolation complexe, au
sens de \cite{BL}. Soit 
$S = \{ z \in \mathbb{C} \, ; \ 0 \leq \Rea (z) \leq 1 \}$.

Rappelons d'abord la d\'{e}finition de l'espace d'interpolation $A_\theta$,
o\`{u} $\theta \in \, ]0, 1[$ \cite[chapitre 4]{BL}. On note $\FdeA$
l'espace des fonctions $F$ \`{a} valeurs dans $A_0 + A_1$, continues
born\'{e}es sur $S$, holomorphes \`{a} l'int\'{e}rieur de $S$, telles que,
pour $j \in \{0, 1\}$, $F(j + i \tau)$ prend ses valeurs dans $A_j$ et 
$\| F(j + i\tau) \|_{A_j} \rightarrow 0$, 
quand $|\tau| \rightarrow +\infty $. On munit $\FdeA$ de la norme

\begin{equation*}
   \|F\|_{\FdeA}
 = \max
  (\sup_{\tau \in \R} \| F(i\tau) \|_{A_0}, \,
   \sup_{\tau \in \R} \| F(1 + i\tau) \|_{A_1}).
\end{equation*}

\noindent L'espace $A_\theta = (A_0, A_1)_\theta
 = \{ F(\theta) ; \ F \in \FdeA \}$ est un Banach \cite[theorem 4.1.2]{BL}
pour la norme d\'{e}finie par

\begin{equation*}
   \|a\|_{A_\theta} 
 = \inf \left\{ \| F \|_{\FdeA} ; \ F(\theta) = a \right\} .
\end{equation*}

Rappelons maintenant la d\'{e}finition de l'espace d'interpolation 
$A^\theta$ \cite[chapitre 4]{BL}. On note $\GdeA$ l'espace des fonctions
$g$ \`{a} valeurs dans $A_0 + A_1$, continues sur $S$, holomorphes \`{a}
l'int\'{e}rieur de $S$, telles que $z \rightarrow
 (1 + |z|)^{-1} \|g(z)\|_{A_0 + A_1}$
est born\'{e}e sur $S$, $g(j + i \tau ) - g(j + i \tau^\prime) \in A_j$
pour tous $\tau, \tau^\prime \in \R$, $j \in \{0, 1\}$, et la
quantit\'{e} suivante est finie:
\begin{eqnarray*}
&&
   \| \dot g \|_{Q \GdeA} \\
&=&
 \max 
 \left[   
  \sup_{\underset{\tau \neq \tau^\prime}{\tau, \tau^\prime \in \R}}
   \left\| \frac{g(i\tau ) - g(i\tau^\prime)}{\tau - \tau^\prime}
   \right\|_{A_0},
   \ 
  \sup_{\underset{\tau \neq \tau^\prime}{\tau, \tau^\prime \in \R}}
   \left\| \frac{g(1 + i\tau) - g(1 + i\tau^\prime)}{\tau - \tau^\prime}
   \right\|_{A_1}
 \right].
\end{eqnarray*}

\noindent Cette quantit\'{e} d\'{e}finit une norme sur l'espace $Q \GdeA$,
quotient de $\GdeA$ par les applications constantes \`{a} valeurs dans 
$A_0 \cap A_1$, et $Q \GdeA$ est complet pour cette norme 
\cite[lemma 4.1.3]{BL}.

On rappelle \cite[p. 89]{BL} que, pour $g \in \GdeA$,

\begin{equation}
      \| g^\prime (z) \|_{A_0 + A_1} 
 \leq \| \dot g \|_{Q \GdeA}, 
 \;\; z \in S. 
 \label{2}
\end{equation}

\noindent C'est une cons\'{e}quence imm\'{e}diate de l'in\'{e}galit\'{e}

\begin{equation*}
      \left\| \frac{g(z + it) - g(z)} {t} \right\|_{A_0 + A_1}
 \leq \| \dot g \|_{Q \GdeA},
      \; z \in S, \; t \in \R^\ast,
\end{equation*}

\noindent qui d\'{e}coule de la d\'{e}finition de $\| \dot g \|_{Q \GdeA}$
et du th\'{e}or\`{e}me des trois droites \cite[lemma 1.1.2]{BL}
appliqu\'{e} aux fonctions $z \rightarrow
 \langle (g(z+it) - g(z)) / t, a^\ast \rangle$, 
$t$ r\'{e}el fix\'{e}, o\`{u} $a^\ast$ parcourt la boule unit\'{e} de 
$A_0^\ast \cap A_1^\ast$.

L'espace $A^\theta = \{ \ g^\prime (\theta);  \ g \in \GdeA \ \}$ 
est un Banach \cite[theorem 4.1.4]{BL} pour la norme d\'{e}finie par

\begin{equation*}
   \|a\|_{A^\theta}
 = \inf \left\{ \| \dot g \|_{Q \GdeA} \, ; \ g^\prime(\theta) = a 
        \right\}.
\end{equation*}

\noindent D'apr\`{e}s (\ref{2}) $\|a\|_{A_0 + A_1} \leq  \|a\|_{A^\theta}$.
La contraction $A^\theta \rightarrow A_0 + A_1$ est injective par
d\'{e}finition de $A^\theta$.

\noindent D'apr\`{e}s \cite{B}, $A_\theta$ s'identifie isom\'{e}triquement 
\`{a} un sous espace de $A^\theta$.

\noindent\ \ \ \ \ D'apr\`{e}s \cite[theorem 4.2.2]{BL}, $A_0 \cap A_1$
est toujours dense dans $A_\theta$, $0 < \theta < 1$. Si $A_0 \cap A_1$ est
dense dans $A_0$ et $A_1$, on a $(A_0 \cap A_1)^\ast = A_0^\ast + A_1^\ast$, 
$A_0^\ast \cap A_1^\ast = (A_0 + A_1)^\ast$ \cite[theorem 2.7.1]{BL}, on
peut appliquer le th\'{e}or\`{e}me d'it\'{e}ration \cite[theorem 4.6.1]{BL}
et $(A_\theta)^\ast = (A_0^\ast, A_1^\ast)^\theta$, $\theta \in \, ]0, 1[$
\cite[theorem 4.5.1]{BL}. On fait cette hypoth\`{e}se dans la suite.

\noindent

\begin{definition}
\cite{DGZ} Un espace de Banach $X$ est localement uniform\'{e}ment convexe,
ce qu'on note LUR (resp. faiblement LUR) si, pour tout $x\in X$ et pour
toute suite $(x_n)_{n \geq 0}$ dans $X$ satisfaisant
\begin{equation*}
 \| x_n \|^2 / 2 + \|x\|^2 / 2 - \| (x_n + x) / 2 \|^2
 \rightarrow_{n \rightarrow \infty} 0,
\end{equation*}

alors $x_n \rightarrow_{n \rightarrow \infty} x$ en norme (resp.
faiblement).
\end{definition}

\section{R\'{e}sultats}

Notons $B_j$ l'adh\'{e}rence de $A_0^\ast \cap A_1^\ast$ dans $A_j^\ast$,
$j = 0, 1$. Il est clair que $B_0 \cap B_1 =  A_0^\ast \cap A_1^\ast$,
isom\'{e}triquement. D'apr\`{e}s \cite[Theorem 4.2.2 b)]{BL} on a
isom\'{e}triquement, pour $\theta \in \, ]0, 1[$,

\begin{equation}
(B_0, B_1)_\theta = (A_0^\ast, A_1^\ast)_\theta \, .
\label{II}
\end{equation}

\noindent Comme $B_0 \cap B_1$ est dense dans $B_j$, le dual de 
$B_\theta = (B_0, B_1)_\theta$ est $(B_0^\ast, B_1^\ast)^\theta$ 
\cite[theorem 4.5.1]{BL} et, d'apr\`{e}s \cite[theorem 2.7.1]{BL},

\begin{equation*}
   B_0^\ast + B_1^\ast
 = (B_0 \cap B_1)^\ast
 = (A_0^\ast \cap A_1^\ast)^\ast
 = (A_0 + A_1)^{\ast \ast}.
\end{equation*}

\noindent En particulier, $A_0 + A_1$ s'identifie isom\'{e}triquement \`{a}
un sous espace ferm\'{e} de $B_0^\ast + B_1^\ast$.

\noindent Soit $i_j : B_j \rightarrow A_j^\ast$ l'application identit\'{e}; 
la restriction de son adjoint $i_j^\ast : A_j \rightarrow B_j^\ast$, 
$j = 0, 1$, est contractante.

\begin{lemma}
\label{R}
Soit $R : Q \Gde{A_0, A_1} \rightarrow Q \Gde{B_0^\ast, B_1^\ast}$ 
l'application qui est d\'{e}finie par 
$g(j + i \, \cdot) \rightarrow i_j^\ast (g(j + i \, \cdot))$, $j = 0, 1$. 
L'application $R$ \ est une contraction et induit une contraction injective
\end{lemma}

\begin{equation*}
 R^\theta : A^\theta \rightarrow (B_0^\ast, B_1^\ast)^\theta,
 \; \theta \in \, ]0, 1[ .
\end{equation*}

\noindent D\'{e}monstration: Il est clair que $R$ est une contraction (non
injective en g\'{e}n\'{e}ral). On identifie $A^\theta$ et 
$(B_0^\ast, B_1^\ast)^\theta$ \`{a} des quotients de $Q \Gde{A_0, A_1}$ 
et $Q \Gde{B_0^\ast, B_1^\ast}$ respectivement. Notant que
$(R( \dot g ))^\prime (\theta) = R^\theta (g^\prime (\theta))$, $R$ induit une
contraction $R^\theta$ sur ces quotients. Notons que, pour 
$a \in A^\theta$, pour $b \in B_0 \cap B_1 = A_0^\ast \cap A_1^\ast
 = (A_0 + A_1)^\ast $ (espace dense dans $B_\theta$),

\begin{equation*}
   \langle R^\theta (a), b \rangle
 = \langle a, b \rangle .
\end{equation*}
Si $R^\theta (a) = 0$, alors $\langle a, b \rangle = 0$ pour tout 
$b \in B_0 \cap B_1 = (A_0 + A_1)^\ast$, d'o\`{u} $a = 0$ dans 
$A_0 + A_1$, et dans $A^\theta$. 
$\blacksquare$

\begin{theorem}
\label{ci}
Soient $(A_0, A_1)$ un couple d'interpolation complexe et $B_j$
l'ad\-h\'{e}rence de $A_0^\ast \cap A_1^\ast$ dans $A_j^\ast$, $j = 0, 1$,
$\beta \in \, ]0, 1[$, $R^\beta$ d\'{e}finie comme ci-dessus. Soit
$Z^\beta$ l'adh\'{e}rence de $R^\beta (A^\beta)$ dans 
$(B_0^\ast, B_1^\ast)^\beta$. Supposons que $Z^\beta$ est un espace
faiblement-LUR. Alors $A^\theta = A_\theta$, pour tout 
$\theta \in \, ]0, 1[$.
\end{theorem}

\noindent La d\'{e}monstration n\'{e}cessite les lemmes suivants.

\begin{lemma}
\label{yu}\textbf{\ }Pour tout $\theta \in \, ]0, 1[$, $R^\theta$ est une
isom\'{e}trie: $A_\theta \rightarrow (B_0^\ast, B_1^\ast)^\theta$.
\end{lemma}

\noindent D\'{e}monstration: Comme $A_\theta$ s'identifie \`{a} un 
sous-espace de $A^\theta$ \cite{B}, $R^\theta$ est contractante: 
$A_\theta = (A_0, A_1)_\theta \rightarrow (B_0^\ast, B_1^\ast)^\theta$
par le lemme \ref{R}. Comme $A_0 \cap A_1$ est dense dans $A_\theta$, il
suffit de montrer que $\|a\|_{A_\theta}
 \leq \| R^\theta (a) \|_{(B_0^\ast, B_1^\ast)^\theta}$ 
lorsque $a \in A_0 \cap A_1$.

\noindent Soit $\varepsilon > 0$; comme 
$(A_\theta)^\ast = (A_0^\ast, A_1^\ast)^\theta$, il existe 
$g \in \Gde{A_0^\ast, A_1^\ast}$ tel que

\begin{equation}
   \|a\|_{A_\theta}
 < \left| \langle a, g^\prime(\theta) \rangle \right| + \varepsilon, 
 \text{ \qquad }
   \| \dot g \|_{Q \Gde{A_0^\ast, A_1^\ast}} \leq 1.
\label{io}
\end{equation}

\noindent Soient

\begin{equation*}
   F_n(z)
 = i n \, [ g(z + i/n) - g(z) ], 
 \quad 
 z \in S
\end{equation*}
et $F_{n, \delta}(z) = e^{\delta z^{2}} F_n(z)$, \ pour $\delta > 0$. Comme 
$|F_n|$ est born\'{e}e sur le bord de $S$, $|F_{n, \delta }|$ tend vers $0$
\`{a} l'infini sur le bord, d'o\`{u}
$F_{n, \delta } \in \Fde{A_0^\ast, A_1^\ast}$. Par d\'{e}finition 
\begin{equation*}
      \| F_{n, \delta }(\theta) \|_{(A_0^\ast, A_1^\ast)_\theta}
 \leq \| F_{n, \delta } \|_{\Fde{A_0^\ast, A_1^\ast }}
 \leq e^{\delta } \sup_{z \in S} \, |F_n(z)|
 \leq e^{\delta } \| \dot g \|_{Q \Gde{A_0^\ast, A_1^\ast}} 
 \leq e^{\delta }.
\end{equation*}

\noindent D'o\`{u}, pour tout $n$, par (\ref{II}),

\begin{equation*}
   \| F_n(\theta) \|_{(B_0, B_1)_\theta}
 = \| e^{-\delta \theta^2} 
       F_{n, \delta }(\theta) \|_{(A_0^\ast, A_1^\ast )_\theta }
 = \lim_{\delta \rightarrow 0}
    \| e^{-\delta \theta^2} 
        F_{n, \delta }(\theta) \|_{(A_0^\ast, A_1^\ast)_\theta }
 \leq 1.
\end{equation*}

\noindent Comme $g$ est holomorphe \`{a} valeurs dans 
$A_0^\ast + A_1^\ast = (A_0 \cap A_1)^\ast$, 
$\langle a, F_n(\theta) \rangle 
 \underset{n \rightarrow \infty}{\rightarrow}
 \left\langle a, g^\prime(\theta) \right\rangle$. Il existe $n_0$ assez
grand tel que, d'apr\`{e}s (\ref{io}),

\begin{eqnarray*}
&&
   -2 \varepsilon + \| a \|_{A_\theta}
 < \left| \langle a, F_{n_0}(\theta) \rangle \right| \\
 &\leq& \| R^\theta (a) \|_{(B_0^\ast, B_1^\ast )^\theta} 
        \, \| F_{n_0}(\theta) \|_{(B_0, B_1)_\theta}
 \leq \| R^\theta (a) \|_{(B_0^\ast, B_1^\ast )^\theta},
\end{eqnarray*}

\noindent d'o\`{u} l'in\'{e}galit\'{e} cherch\'{e}e lorsque 
$\varepsilon \rightarrow 0$. 
$\blacksquare$

\begin{lemma}
\label{bn} Soient $g \in \GdeA$, $\theta \in \, ]0, 1[$. L'application:
$\tau \rightarrow R^\theta (g^\prime (\theta + i\tau ))$ est born\'{e}e de
$\R$ dans $(B_0^\ast, B_1^\ast)^\theta$. Pour tout 
$c \in (B_0^\ast, B_1^\ast )^\theta$, l'application: $\tau \rightarrow
 \left\| c + R^\theta (g^\prime (\theta + i\tau ))
 \right\|_{(B_0^\ast, B_1^\ast )^\theta}$ est s.c.i. sur $\R$.
\end{lemma}

\noindent D\'{e}monstration: Par d\'{e}finition de $A^\theta$, 
$g^\prime(\theta) \in A^\theta$ ; par le lemme \ref{R}

\begin{equation*}
      \| R^\theta (g^\prime(\theta)) \|_{(B_0^\ast, B_1^\ast )^\theta}
 \leq \| g^\prime (\theta) \|_{A^\theta }
 \leq \| \dot g \|_{Q \GdeA}.
\end{equation*}

\noindent La fonction $g_{i\tau}$ d\'{e}finie par 
$g_{i\tau }(z) = g(z + it)$, $z \in S$, $\tau \in \R$, v\'{e}rifie 
$\| \dot g_{i\tau} \|_{Q \GdeA} = \| \dot g \|_{Q \GdeA}$, 
donc $\left\| R^\theta (g_{i\tau}^\prime (\theta))
      \right\|_{(B_0^\ast, B_1^\ast)^{\theta}}
      \leq \| \dot g \|_{Q \GdeA}$.

\noindent D'apr\`{e}s (\ref{II}), et comme $B_0 \cap B_1
 = A_0^\ast \cap A_1^\ast$ est dense dans $B_\theta$, on a

\begin{eqnarray*}
&&
  \| c + R^\theta 
          (g^\prime(\theta + i\tau )) \|_{(B_0^\ast, B_1^\ast )^\theta} \\
&=&
  \sup 
   \left\{ 
    \left| \langle b, c + R^\theta (g^\prime (\theta + i\tau ))
           \rangle 
    \right| ;
    \text{ }
    \| b \|_{(B_0, B_1)_\theta } \leq 1
   \right\} \\
&=&
  \sup
   \left\{ 
    \left| \langle a^\ast, c+g^\prime (\theta + i\tau)
           \rangle 
    \right| ;
    \text{ }
    a^\ast \in A_0^\ast \cap A_1^\ast, 
    \quad 
    \| a^\ast \|_{(A_0^\ast, A_1^\ast )_\theta} \leq 1
   \right\}.
\end{eqnarray*}
Comme $g$ est holomorphe \`{a} valeurs dans $A_0 + A_1$, pour tout 
$a^\ast \in A_0^\ast \cap A_1^\ast = (A_0 + A_1)^\ast$, les
applications $\tau \rightarrow
 \left| \langle a^\ast, c + g^\prime (\theta + i\tau) \rangle \right|$ sont
continues sur $\R$. Leur supremum est donc s.c.i.. 
$\blacksquare$

\begin{lemma}
\label{k}
Soient $\overline{C} = (C_0, C_1)$ un couple d'interpolation, 
$\beta \in \, ]0, 1[$, $Z^\beta$ un sous-espace ferm\'{e}
faiblement-LUR de $C^\beta$, $g \in \Gde{\overline{C}}$. On suppose que
l'application $\phi_\beta : \tau \in \R 
 \rightarrow g^\prime (\beta + i\tau)$ est born\'{e}e \`{a} valeurs dans
$Z^\beta$ et que l'application $\| c + \phi_\beta \|_{C^\beta}$ est s.c.i.,
pour tout $c \in Z^\beta$ fix\'{e}. Alors $\phi_\beta$ est p.s. \'{e}gale
\`{a} une fonction fortement mesurable: 
$\R \rightarrow Z^\beta \subset C^\beta$.
\end{lemma}

\noindent Preuve: Comme $s \rightarrow \| \phi_\beta (\tau )
 + \phi_\beta (s) \|_{Z^\beta}$ est s.c.i. sur $\R$, si 
$\tau_n \rightarrow \tau $,

\begin{eqnarray*}
&&
      0
 \leq \overline{\lim}_{n \rightarrow +\infty } E_n \\
&=&
  \overline{\lim}
   \left\{ 2 \| \phi_\beta (\tau) \|_{Z^\beta}^2
    + 2 \| \phi_\beta (\tau_n) \|_{Z^\beta}^2 
    - \| \phi_\beta (\tau) + \phi_\beta (\tau_n) \|_{Z^\beta}^2 
   \right\} \\
&
 \leq & 2 \, \| \phi_\beta (\tau) \|_{Z^\beta }^2
  + 2 \, \overline{\lim} \| \phi_\beta (\tau_n) \|_{Z^\beta}^2
  - \underline{\lim} 
     \| \phi_\beta (\tau) + \phi_\beta (\tau_n) \|_{Z^\beta}^2 \\
&
 \leq & 2 \, \| \phi_\beta (\tau ) \|_{Z^\beta}^2
  + 2 \, \overline{\lim } \| \phi_\beta (\tau_n) \|_{Z^\beta}^2
  - 4 \, \| \phi_\beta (\tau) \|_{Z^\beta}^2 \\
&=&
   2 \, \overline{\lim } \| \phi_\beta (\tau_n) \|_{Z^\beta}^2
  -2 \, \| \phi_\beta (\tau) \|_{Z^\beta}^2.
\end{eqnarray*}

\noindent Comme $\| \phi_\beta \|_{Z^\beta}$ est mesurable born\'{e}e, pour
tout $N$ et $\varepsilon > 0$, il existe, d'apr\`{e}s le th\'{e}or\`{e}me
de Lusin, un compact $K_{N, \varepsilon} \subset [-N, N]$, de mesure 
$> 2N - \varepsilon$, sur lequel $\| \phi_\beta \|_{Z^\beta}$ est continue.
Soit $(\tau_n)_{n \geq 0}$ une suite dans $K_{N, \varepsilon}$ convergeant 
vers $\tau$. D'apr\`{e}s ce qui pr\'{e}c\`{e}de 
$\overline{\lim }_{n \rightarrow +\infty } E_n = 0$. Comme $E_n \geq 0$,
$E_n \rightarrow_{n \rightarrow \infty} 0$. Par d\'{e}finition 
de la propri\'{e}t\'{e} faiblement-LUR de $Z^\beta$, cela entra\^{\i}ne que
$\phi_\beta (\tau_n) \rightarrow \phi_\beta (\tau)$ faiblement dans
$Z^\beta$, c\`{a}d $\phi_\beta$ est faiblement continue sur 
$K_{N, \varepsilon}$. Soit $Y$ le sous espace ferm\'{e} de $Z^\beta$
engendr\'{e} par $\phi_\beta (K_{N, \varepsilon})$. Alors $Y$ est
s\'{e}parable: sinon, \'{e}tant donn\'{e}e une suite $(s_n)_{n \geq 1}$
dense dans $K_{N, \varepsilon}$, il existe, d'apr\`{e}s le th\'{e}or\`{e}me
de Hahn-Banach, $z \in Y^\ast $, non nul, tel que 
$(\phi_\beta (s_n), z) = 0$ pour tout $n$; par continuit\'{e} 
$s \rightarrow (\phi_\beta (s), z)$ est nulle sur $K_{N, \varepsilon}$,
d'o\`{u} $z = 0$ et la contradiction. Par le th\'{e}or\`{e}me de Pettis
\cite[theorem II 2]{DU}, $\phi_\beta$ est fortement mesurable:
$K_{N, \varepsilon } \rightarrow Y \subset Z^\beta$. Cela montre le
r\'{e}sultat annonc\'{e}. 
$\blacksquare$

\begin{lemma}
\label{mesurable}
Soient $g \in \GdeA$, $\phi_\theta (t) = g^\prime (\theta + it)$, 
$t \in \R$.

i) Si $\phi_\theta$ est \`{a} valeurs dans un sous espace ferm\'{e} 
s\'{e}parable $Z$ de $A_\theta$, elle est fortement mesurable: 
$\R \rightarrow A_\theta$.

ii) Si $\phi_\theta$ est \`{a} valeurs dans un sous espace ferm\'{e} 
s\'{e}parable $Z$ de $A^\theta$, elle est fortement mesurable: 
$\R \rightarrow A^\theta$.
\end{lemma}

\noindent Dans la suite on utilise seulement i), dans la preuve du lemme
\ref {lem} d). On donne deux preuves de i) (noter que ii) implique i)).

\noindent Preuve: i) D'apr\`{e}s le lemme \ref{yu}, $Z$ est un sous espace
ferm\'{e} de $(B_0^\ast, B_1^\ast )^\theta$ et, d'apr\`{e}s le
lemme \ref{bn}, l'application $t \rightarrow \| \phi_\theta (t) - c \|_Z$
est s.c.i. pour tout $c \in Z$. L'image r\'{e}ciproque par $\phi_\theta$ de
toute boule ouverte de $Z$ est donc un bor\'{e}lien. Comme $Z$ est
s\'{e}parable, tout ouvert de $Z$ est r\'{e}union d\'{e}nombrable de
boules, donc $\phi_\theta$ est bien mesurable \`{a} valeurs dans $Z$. 
$\blacksquare$

ii) Soient $J$ l'injection canonique: $Z \rightarrow A_0 + A_1$, et $Y$
l'adh\'{e}rence de $J(Z)$ dans $A_0 + A_1$. Comme $Z$ et $Y$ sont des
espaces polonais, comme $J$ est continue, $J^{-1}$ est bor\'{e}lienne: 
$J(Z) \rightarrow Z$, voir par exemple \cite{A}. Comme 
$J \circ \phi_\theta : \R \rightarrow A_0 + A_1$ est continue et \`{a}
valeurs dans $J(Z)$, comme 
$\phi_\theta = J^{-1} \circ (J \circ \phi_\theta)$, alors $\phi_\theta$ est
bor\'{e}lienne: $\R \rightarrow Z$. 
$\blacksquare$

\begin{lemma}
\label{lem} Soient $g \in \GdeA$, $\beta \in \, ]0, 1[$ 
et $\phi_\beta (\cdot) = g^\prime (\beta + i \, \cdot)$.

a) On suppose que $R^\beta \circ \phi_\beta$ est p.s. \'{e}gale \`{a}
une fonction fortement mesurable: 
$\R \rightarrow (B_0^\ast, B_1^\ast)^\beta$. Alors $\phi_\beta$ est p.s.
\`{a} valeurs dans $A_\beta$.

On suppose d\'{e}sormais que $\phi_\beta$ est p.s. \'{e}gale \`{a} une
fonction fortement mesurable: $\R \rightarrow A_\beta$. Alors

b) pour $\theta \neq \beta$, $g^\prime (\theta) \in A_\theta$.

c) pour tout $\theta \neq \beta$, $\phi_\theta$ est \`{a} valeurs dans
un sous espace s\'{e}parable de $A_\theta$.

d) $g^\prime (\beta ) \in A_\beta$.
\end{lemma}

\noindent On a not\'{e} $R^\beta \circ \phi_\beta$ la fonction: 
$t \rightarrow R^\beta (g_{it}^\prime (\beta))$.

\noindent Preuve: a) \emph{\'{e}tape 1: } Comme $g$ est holomorphe \`{a}
l'int\'{e}rieur de $S$, pour tous $t \in \R$, $h > 0$, $\theta \in \, ]0, 1[$,
on a, dans $A_0 + A_1$,

\begin{equation}
   g(\theta + i(t+h)) - g(\theta + it)
 = \int_t^{t+h} g^\prime (\theta + i\tau) \, d\tau  
 \label{4bis}
\end{equation}

\noindent Posons

\begin{equation*}
   g_1
 = g - g(0) - \alpha_0
\end{equation*}
o\`{u} $g(1) - g(0) = \alpha_0 + \alpha_1$ ($\alpha_j \in A_j$, $j = 0, 1$),
avec

\begin{equation*}
   \| g(1) - g(0) \|_{A_0 + A_1}
 = \| \alpha_0 \|_{A_0} + \| \alpha_1 \|_{A_1}.
\end{equation*}

\noindent D'apr\`{e}s l'in\'{e}galit\'{e} des accroissements finis et 
(\ref{2})

\begin{equation*}
      \| g(1) - g(0) \|_{A_0 + A_1}
 \leq \| \dot g \|_{Q \GdeA}.
\end{equation*}
Alors $g_1 : S \rightarrow A_0 + A_1$ est continue sur $S$ et holomorphe
\`{a} l'int\'{e}rieur de $S$. Comme $g \in \GdeA$, pour tout $\tau \in \R$
et $j\in \{0, 1\}$, on a

\begin{equation*}
      \| g_1 (j + i\tau) \|_{A_j}
 \leq \| g(j + i\tau) - g(j) \|_{A_j} + \| \alpha_j \|_{A_j }
 \leq (1 + |\tau| ) \| \dot g \|_{Q \GdeA}.
\end{equation*}

\noindent L'application $z \rightarrow G_\varepsilon (z)
 = e^{\varepsilon z^{2}} g_1 (z)$ est donc dans $\FdeA$ pour tout 
$\varepsilon > 0$. En particulier, pour tout $t \in \R$, 
$G_\varepsilon (\theta + it) \in A_\theta$, donc 
$g_1 (\theta + it) \in A_\theta$. D'o\`{u}

\begin{equation*}
   g_1 (\theta + i(t+h)) - g_1 (\theta + it)
 = g(\theta + i(t+h)) - g(\theta + it) \in A_\theta.
\end{equation*}

\noindent Alors, d'apr\`{e}s (\ref{4bis}), 
$\int_t^{t+h} g^\prime (\theta + i\tau) \, d\tau$ est dans $A_\theta$, pour
$t$ et $h$ r\'{e}els.

\emph{\'{e}tape 2: } Par hypoth\`{e}se $R^\beta \circ \phi_\beta$ est
p.s. \'{e}gale \`{a} une fonction fortement mesurable: 
$\R \rightarrow Z^\beta$, o\`{u} $Z^\beta$ est l'adh\'{e}rence de $A^\beta$
dans $(B_0^\ast,  B_1^\ast )^\beta$. Le th\'{e}or\`{e}me de
diff\'{e}rentiabilit\'{e} de Lebesgue \cite[chap. II theorem 9 p 48]{DU}
entra\^{\i}nent que, p.s., on a dans $Z^\beta$ l'\'{e}galit\'{e}

\begin{equation}
   i R^\beta \circ \phi_\beta (it)
 = \lim_{h \rightarrow 0}
    \frac{1}{h}
    \int_t^{t+h} R^\beta \circ \phi_\beta (i\tau ) \, d\tau
 = \lim_{h \rightarrow 0} 
    R^\beta \bigl( \frac{1}{h} \int_t^{t+h} g^\prime (\beta + i\tau ) \, d\tau
            \bigr),  
 \label{Leb}
\end{equation}

\noindent o\`{u} $h$ est r\'{e}el. D'apr\`{e}s la fin de l'\'{e}tape~1
appliqu\'{e}e en $\beta$ et le lemme \ref{yu}, cette limite dans $Z^\beta$
est en fait une limite dans $A_\beta$, c\`{a}d p.s. 
$g^\prime (\beta + i \, \cdot) \in A_\beta$.

b) On suppose d'abord $\theta > \beta$.

\noindent \emph{\'{e}tape~1: } Soit

\begin{equation*}
 V(z) = g_1 (\beta + (1-\beta)z), \; z \in S.
\end{equation*}

\noindent Cette fonction \`{a} valeurs dans $A_0 + A_1$ est holomorphe
\`{a} l'int\'{e}rieur de $S$ et continue sur $S$, donc s'exprime \`{a}
l'aide de la mesure harmonique sur le bord de $S$. Pour v\'{e}rifier que
$V$, vue comme fonction \`{a} valeurs dans $A_\beta + A_1$, est holomorphe
\`{a} l'int\'{e}rieur de $S$ et continue sur $S$, il suffira donc de voir
que $V$ est continue sur l'axe imaginaire, \`{a} valeurs dans $A_\beta$.

On va montrer que $V \in \Gde{A_\beta, A_1}$ avec une norme 
$\leq (1-\beta ) \| \dot g \|_{Q \GdeA}$. L'in\'{e}galit\'{e}
correspondante sur la droite $\Rea z = 1$ est \'{e}vidente. Pour la
v\'{e}rifier sur l'axe imaginaire, posons, pour $\tau, \tau^{\prime}$ 
r\'{e}els fix\'{e}s,

\begin{equation*}
   F_{\tau, \tau^\prime}(\xi)
 = \frac{g(\xi + i(1-\beta )\tau) - g(\xi + i(1-\beta )\tau^\prime)}
    {\tau -\tau^\prime }, \; \xi \in S,
\end{equation*}

\noindent d'o\`{u} $F_{\tau, \tau^\prime }(\beta)
 = (V(i\tau) - V(i\tau^\prime )) / (\tau - \tau^\prime)$, 
et $F_{\tau, \tau^\prime}(1)
 = (V(1 + i\tau ) - V(1 + i\tau^\prime )) / (\tau - \tau^\prime)$. Pour
tout $t \in \R$, on a

\begin{equation*}
      \| F_{\tau, \tau^\prime} (j + it) \|_{A_j}
 \leq (1 - \beta) \| \dot g \|_{Q \GdeA},
  \ j \in \{0, 1\}.
\end{equation*}

\noindent Comme dans l'\'{e}tape~1 de a), pour tout $\varepsilon > 0$,
l'application $\xi \rightarrow H_{\varepsilon, \tau, \tau^\prime}(\xi)
 = e^{\varepsilon \xi^{2}} F_{\tau, \tau^\prime}(\xi)$ v\'{e}rifie

\begin{equation*}
      \| H_{\varepsilon, \tau, \tau^\prime} \|_{\FdeA}
 \leq e^{\varepsilon }(1 - \beta) \| \dot g \|_{Q \GdeA},
\end{equation*}

\noindent d'o\`{u}

\begin{equation*}
      \| F_{\tau, \tau^\prime }(\beta) \|_{A_\beta}
 \leq (1 - \beta ) \| \dot g \|_{Q \GdeA}.
\end{equation*}

\noindent On a donc, pour tous $\tau, \tau^\prime$ r\'{e}els,

\begin{equation*}
      \| V(i\tau ) - V(i\tau^\prime) \|_{A_\beta}
 \leq |\tau - \tau^\prime| (1 - \beta) \| \dot g \|_{Q \GdeA},
\end{equation*}

\noindent ce qui prouve la continuit\'{e} de $V$ sur l'axe imaginaire, \`{a}
valeurs dans $A_\beta $, et l'assertion annonc\'{e}e.

\emph{\'{e}tape 2: } d'apr\`{e}s la preuve de a), pour $h$ r\'{e}el, p.s.

\begin{equation*}
   \lim_{h \rightarrow 0} (V(i(\tau + h)) - V(i\tau )) / h
 = (1 - \beta) g^\prime (\beta + (1-\beta ) i\tau)
   \quad \text{dans} \quad A_\beta .
\end{equation*}

\noindent D'apr\`{e}s \cite[lemma 4.3.3]{BL}, on a alors

\begin{equation*}
 V^\prime (\eta ) \in (A_\beta, A_1)_\eta,
  \eta \in \, ]0, 1[ .
\end{equation*}

\emph{\'{e}tape 3:} Choisissons $\eta$ tel que $\theta
 = (1 - \eta) \beta + \eta $. D'apr\`{e}s le th\'{e}or\`{e}me de
r\'{e}it\'{e}ration \cite[theorem 4.6.1]{BL}, 
$(A_\beta, A_1)_\eta = A_\theta$, donc

\begin{equation*}
 V^\prime (\eta ) = (1 - \beta) g^\prime (\theta) \in A_\theta,
\end{equation*}

\noindent ce qui ach\`{e}ve la preuve lorsque $\beta < \theta$.

Si $0 < \theta < \beta $ le raisonnement est analogue, en rempla\c{c}ant $V$
par $W(z) = g_1 (\beta z) \in \Gde{A_0, A_\beta}$, telle que 
$\lim_{h \rightarrow 0} (W(1 + i(\tau + h)) - W(1 + i\tau )) / h$ existe
dans $A_\beta$, pour presque tout $\tau$, avec $h$ r\'{e}el.

c) Soit $A_0^\prime \subset A_0$ \ le sous espace ferm\'{e} s\'{e}parable 
engendr\'{e} par $\{g_1 (it), \, t \in \R \}$. Comme $g_1$ est continue sur
$S$, $A_0^\prime$ est s\'{e}parable, ainsi que $(A_0^\prime, A_1)_\beta$ et
son adh\'{e}rence $Y$ dans $A_\beta$. Par l'\'{e}tape~2 de b) appliqu\'{e}e
au couple $(A_0^\prime, A_1)$, $g^\prime (\beta + it)$ est p.s. dans
$(A_0^\prime, A_1)_\beta$, donc p.s. dans $Y$, ce qui r\`{e}gle le cas
$\theta  = \beta$.

Pour le cas $\beta < \theta$, rempla\c{c}ons la fonction $V$ de 
l'\'{e}tape~1 de b) par $V_t(z) = V(z + it)$, avec $t$ fix\'{e} r\'{e}el.
Comme en b), $V_t \in \Gde{Y, A_1}$, 
$V_t^\prime (\eta ) \in (Y, A_1)_\eta$, $\eta \in \, ]0, 1[ $ et 
$(Y, A_1)_\eta$ est s\'{e}parable. Soit $\eta$ d\'{e}fini comme dans
l'\'{e}tape 3 de b). Comme ci-dessus, 
$V_t^\prime (\eta) = (1 - \beta) g^\prime (\theta + i(1-\beta )t)$. Soit 
$Z_\theta$ l'adh\'{e}rence de $(Y, A_1)_\eta$ dans $(A_\beta, A_1)_\eta
 = A_\theta$; $Z_\theta$ est donc s\'{e}parable et 
$\phi_\theta = g^\prime (\theta + i.)$ est \`{a} valeurs dans $Z_\theta$.

On raisonne de fa\c{c}on analogue si $0 < \theta < \beta$ en consid\'{e}rant 
$W_t(z) = W(z + it)$: $W_t$ est dans $\Gde{A_0^\prime, Y}$.

d) Soit $\theta > \beta$. Par c) et le lemme \ref{mesurable} i), 
$\phi_\theta$ est fortement mesurable \`{a} valeurs dans $A_\theta$. Alors
b) appliqu\'{e} en \'{e}changeant les r\^{o}les de $\beta $ et $\theta$
donne $g^\prime (\beta) \in A_\beta$. 
$\blacksquare$

\smallskip

\noindent \emph{D\'{e}monstration du th\'{e}or\`{e}me 1:} Soient 
$a \in A^\beta$ et $g \in \GdeA$ tel que $a = g^\prime (\beta)$. 
D'apr\`{e}s le lemme \ref{bn}, l'application $R^\beta \circ \phi_\beta:
 \tau \in \R \rightarrow R^\beta (g^\prime (\beta + i\tau))$ 
est \`{a} valeurs dans $Z^\beta \subset (B_0^\ast, B_1^\ast)^\beta$ 
et v\'{e}rifie les hypoth\`{e}ses du lemme \ref{k} pour 
$C^\beta = (B_0^\ast, B_1^\ast )^\beta$. Gr\^{a}ce \`{a} l'hypoth\`{e}se 
sur $Z^\beta$, on peut appliquer le lemme \ref{k}, donc 
$R^\beta \circ \phi_\beta$ est p.s. \'{e}gale \`{a} une fonction fortement
mesurable \`{a} valeurs dans $(B_0^\ast, B_1^\ast)^\beta$, et 
$\phi_\beta$ est p.s. \'{e}gale \`{a} une fonction fortement mesurable 
\`{a} valeurs dans $A_\beta$ par le lemme \ref{lem} a). D'apr\`{e}s le
lemme \ref{lem} b) $g^\prime (\theta) \in A_\theta$ pour tout 
$\theta \neq \beta$. Il en r\'{e}sulte que $A^\theta = A_\theta$, pour tout 
$\theta \neq \beta$. Enfin par le lemme \ref{lem} d) 
$g^\prime (\beta) = a \in A_\beta$, d'o\`{u} $A^\beta = A_\beta$. 
$\blacksquare$

\begin{proposition}
Soient $A_0, A_1$ deux espaces de Banach tels que $A_0$ s'injecte
continuement dans $A_1$, et $\beta \in \, ]0, 1[$. Si $A_\beta$ a la
propri\'{e}t\'{e} de Radon-Nikodym analytique (d\'{e}finie par exemple dans
\cite{DU}) pour un $\beta \in \, ]0, 1[$, alors $A_\theta = A^\theta$ pour
tout $\theta \in \, ]0, 1[$.
\end{proposition}

\noindent Pour $\beta = 1$ ce r\'{e}sultat est \cite[Proposition 3.1]{HP};
appliqu\'{e} au couple $(A_0, A_\beta)$, il donne la conclusion pour 
$\theta \in \, ]0, \beta[$.

\noindent Preuve: D'apr\`{e}s le lemme \ref{lem} b), d), il suffit de
montrer que pour toute $g \in \GdeA$, $\phi_\beta$ est p.s. mesurable \`{a}
valeurs dans $A_\beta$.

On a mentionn\'{e} dans la preuve du lemme \ref{lem} b) que la fonction 
$z \rightarrow W(z) = g_1 (\beta z)$ est dans $\Gde{A_0, A_\beta}$.
\`{A} l'int\'{e}rieur de $S$, $W^\prime$ est donc holomorphe \`{a} valeurs
dans $A_0 + A_\beta = A_\beta$; par (\ref{2}) elle est born\'{e}e. Comme 
$A_\beta$ poss\`{e}de la propri\'{e}t\'{e} de Radon-Nikodym analytique, 
$W^\prime$ admet p.s. des limites non tangentielles au bord de $S$. Soit 
$\psi$ la limite p.s. (dans $A_\beta$) de $W^\prime$ sur la droite 
$\Rea z = 1$; $\psi $ est donc p.s. mesurable \`{a} valeurs dans $A_\beta$.
Comme $g^\prime$ est continue (\`{a} valeurs dans $A_0 + A_1 = A_1$) sur
$S$, $\psi$ coincide p.s. avec la fonction 
$t \rightarrow \beta g^\prime (\beta + i \beta t)$, ce qui ach\`{e}ve la
preuve. 
$\blacksquare$

\begin{corollary}
Si $A_0$ s'injecte continuement dans $A_1$ avec image dense, si $A_\beta$
est un treillis de Banach, et si $(A_0^\ast, A_1^\ast)^\beta$ admet une
norme \'{e}quivalente LUR pour un $\beta \in \, ]0, 1[$, alors
$(A_0^\ast, A_1^\ast)_\theta = (A_0^\ast, A_1^\ast )^\theta$ pour tout
$\theta \in \, ]0, 1[$.
\end{corollary}

\noindent Preuve: Comme $\ell^\infty$ n'admet aucune norme \'{e}quivalente 
LUR \cite[Chap. II, theorem 7.10]{DGZ}, 
$(A_\beta)^\ast = (A_0^\ast, A_1^\ast )^\beta$ ne contient pas $\ell^\infty$
isomorphiquement. Alors, d'apr\`{e}s un r\'{e}sultat de Bessaga-Pelczy\'{n}ski 
\cite[Corollary I 6]{DU}, l'espace $(A_0^\ast, A_1^\ast)^\beta$ ne
contient pas $c_{0}$ isomorphiquement; comme c'est un treillis, il
poss\`{e}de la propri\'{e}t\'{e} de Radon-Nikodym analytique \cite{E}. Son
sous-espace isom\'{e}trique $(A_0^\ast, A_1^\ast )_\beta$ conserve cette
propri\'{e}t\'{e}. \ La proposition pr\'{e}c\'{e}dente appliqu\'{e}e \`{a}
$A_1^\ast, A_0^\ast$ ach\`{e}ve la preuve.
$\blacksquare$

\textbf{Remerciement: }\textit{Je remercie chaleureusement F. Lust-Piquard
pour ses conseils lors de la r\'{e}daction de ce travail. }

\end{document}